\documentclass[12pt]{article}
\usepackage{amsmath}
\usepackage{amsfonts}
\usepackage{mitpress}
\usepackage{graphicx}

\begin{document}

\title{Galerkin method for the numerical solution of the Burgers' equation
by using exponential B-splines}
\author{M. Zorsahin Gorgulu, I. Dag and D. Irk \\
Department of Mathematics-Computer, Eskisehir Osmangazi University, 26480,
Eskisehir, Turkey.}
\maketitle

\begin{abstract}
In this paper, the exponential B-spline Galerkin \ method is set up for
getting the numerical solution of the Burgers' equation. Two numerical
examples\ related to shock wave propogation and travelling wave are studied
to illustrate the accuracy and the efficiency of the method. Obtained
results are compared with some early studies.

\end{abstract}

\section{Introduction}

Burgers' equation in which convection and diffusion play an important role
arises in applications such as meteorology, turbulent flows, modelling of
the shallow water. Burgers' equation is considered to be useful model for
many physically problems. Thus It is often studied for testing of both real
life problems and computational techniques. Not only does exact solutions of
nonlinear convective problem develops discontinuities in finite time, and
might display complex structure near discontinuities. Efficient and accurate
methods are in need to be tackled the complex solutions of the Burgers'
equation, as expected, many numerical researches are strived to beat
difficulties. Though analytical solutions of the Burgers' equation exist for
simple initial condition ,\ the numerical techniques are of interest to meet
requirement of the wide range of solutions of the Burgers' equations. Some
variants of the Spline methods have set up to find the numerical solutions
of the Burgers' equations such as Galerkin finite element method \cite%
{da,dag3,ozis,bs1,chap}, least square method \cite{sk}, collocation method
\cite{ali,dag,ra,dag4,bs2,bs3,mittal}, finite difference method \cite{f3,f1}%
, differential quadrature method \cite{korkmaz,arora}, method based on the
cubic B-spline quasi interpolant \cite{zhu,jiang}, Taylor--Galerkin and
Taylor-collocation methods \cite{dag2}, etc.

Finite element methods are mainly used methods to have good functional
solutions of the differential equations. The accuracy of the finite element
solutions are increased by the selection of suitable basis\ function for the
approximate function over the finite intervals. We will form the combination
of the exponential B-spline as an approximate function for the finite
element method to get the solution of the Burgers' equation. The exponential
B-splines are suggested to interpolate data and function exhibiting sharp
variations \cite{mc1,mc2,mc3,mc4} since polynomial B-splines\ based
interpolation cause unwanted osculation for interpolation. Some solutions of
the Burgers' equation show sharpness. Thus we will construct the finite
element method together with the exponential B-splines\ to have solutions of
the Burgers' equation. Over the finite element, the Galerkin method will be
employed to determine the unknown of the approximate solution. A few
exponential B-spline numerical methods have suggested for some partial
differential equations: Exponential B-spline collocation method is build up
to compute numerical solution of the convection diffusion equation \cite%
{mohammadi}, the Korteweg-de Vries equation \cite{oz} and generalized
Burgers-Fisher equation \cite{bf}.

In this study, we will consider the Burgers' equation
\begin{equation}
u_{t}+uu_{x}-\nu u_{xx}=0  \label{1}
\end{equation}%
where subscripts $x$ and $t$ are space and time parameters, respectively and
$\nu $ is the viscosity coefficient. Boundary conditions of the Eq.(\ref{1})
are chosen from

\begin{equation}
\begin{tabular}{l}
$u\left( a,t\right) =\beta _{1},$ $u\left( b,t\right) =\beta _{2},$ \\
$u_{x}\left( a,t\right) =0,\text{ }u_{x}\left( b,t\right) =0,$ $t\in \left(
0,T\right] $%
\end{tabular}
\label{2}
\end{equation}%
and initial condition is%
\begin{equation}
u\left( x,0\right) =f\left( x\right) ,x\in \left[ a,b\right] .  \label{2a}
\end{equation}%
$f\left( x\right) $ and $\beta _{1},\beta _{2}$ constants are described in
the computational section.

\section{Exponential B-splines Galerkin Finite Element Solution}

Divide spatial interval $[a,b]$ in $N$ subintervals of length $h=\dfrac{b-a}{%
N}$ and $x_{i}=x_{0}+ih$ at the knots $x_{i},$ $i=0,..,N$ and time interval $%
[0,T]$ in $M$ interval of length $\Delta t.$

Let $\phi _{i}\left( x\right) $ be the exponential B-splines defined at the
knots $x_{i},$ $i=0,\ldots ,N,$ together with fictitious knots $x_{i},$ $%
i=-3,-2,-1,N+1,N+2,N+3$ outside the interval $[a,b].$ The $\phi _{i}\left(
x\right) ,$ $i=-1,\ldots ,N+1$ can be defined as

\begin{equation}
\phi _{i}\left( x\right) =\left \{
\begin{array}{lll}
b_{2}\left[ \left( x_{i-2}-x\right) -\dfrac{1}{p}\left( \sinh \left( p\left(
x_{i-2}-x\right) \right) \right) \right] & \text{ \ } & \text{if }x\in \left[
x_{i-2},x_{i-1}\right] ; \\
a_{1}+b_{1}\left( x_{i}-x\right) +c_{1}e^{p\left( x_{i}-x\right)
}+d_{1}e^{-p\left( x_{i}-x\right) } & \text{ } & \text{if }x\in \left[
x_{i-1},x_{i}\right] ; \\
a_{1}+b_{1}\left( x-x_{i}\right) +c_{1}e^{p\left( x-x_{i}\right)
}+d_{1}e^{-p\left( x-x_{i}\right) } & \text{ } & \text{if }x\in \left[
x_{i},x_{i+1}\right] ; \\
b_{2}\left[ \left( x-x_{i+2}\right) -\dfrac{1}{p}\left( \sinh \left( p\left(
x-x_{i+2}\right) \right) \right) \right] & \text{ } & \text{if }x\in \left[
x_{i+1},x_{i+2}\right] ; \\
0 & \text{ } & \text{otherwise}%
\end{array}%
\right.  \label{3}
\end{equation}%
where%
\begin{equation*}
\begin{array}{l}
p=\underset{0\leq i\leq N}{\max }p_{i},\text{ }s=\sinh \left( ph\right) ,%
\text{ }c=\cosh \left( ph\right) , \\
a_{1}=\dfrac{phc}{phc-s},\text{ }b_{1}=\dfrac{p}{2}\left[ \dfrac{c\left(
c-1\right) +s^{2}}{\left( phc-s\right) \left( 1-c\right) }\right] ,\text{ }%
b_{2}=\dfrac{p}{2\left( phc-s\right) }, \\
c_{1}=\dfrac{1}{4}\left[ \dfrac{e^{-ph}\left( 1-c\right) +s\left(
e^{-ph}-1\right) }{\left( phc-s\right) \left( 1-c\right) }\right] ,\text{ }%
d_{1}=\dfrac{1}{4}\left[ \dfrac{e^{ph}\left( c-1\right) +s\left(
e^{ph}-1\right) }{\left( phc-s\right) \left( 1-c\right) }\right] .%
\end{array}%
\end{equation*}

Each basis function $\phi _{i}\left( x\right) $ is twice continuously
differentiable. Table 1 shows the values of $\phi _{i}\left( x\right) ,$ $%
\phi _{i}^{\prime }\left( x\right) $ and $\phi _{i}^{\prime \prime }\left(
x\right) $ at the knots $x_{i}$:%
\begin{equation*}
\begin{tabular}{c|ccccc}
\multicolumn{6}{l}{Table 1: Exponential B-spline values} \\ \hline \hline
& $x_{i-2}$ & $x_{i-1}$ & $x_{i}$ & $x_{i+1}$ & $x_{i+2}$ \\ \hline
$\phi _{i}\left( x\right) $ & $0$ & $\frac{s-ph}{2\left( phc-s\right) }$ & $%
1 $ & $\frac{s-ph}{2\left( phc-s\right) }$ & $0$ \\
$\phi _{i}^{\prime }\left( x\right) $ & $0$ & $\frac{p\left( c-1\right) }{%
2\left( phc-s\right) }$ & $0$ & $\frac{p\left( 1-c\right) }{2\left(
phc-s\right) }$ & $0$ \\
$\phi _{i}^{\prime \prime }\left( x\right) $ & $0$ & $\frac{p^{2}s}{2\left(
phc-s\right) }$ & $\frac{-p^{2}s}{phc-s}$ & $\frac{p^{2}s}{2\left(
phc-s\right) }$ & $0$ \\ \hline \hline
\end{tabular}%
\end{equation*}

The $\phi _{i}\left( x\right) ,i=-1,\ldots ,N+1$ form a basis for functions
defined on the interval $[a,b]$. The Galerkin method consists of seeking
approximate solution in the following form:%
\begin{equation}
u\left( x,t\right) \approx U_{N}\left( x,t\right) =\overset{N+1}{\underset{%
i=-1}{\sum }}\phi _{i}\left( x\right) \delta _{i}\left( t\right)  \label{4}
\end{equation}%
where $\delta _{i}\left( t\right) $ are time dependent unknown to be
determined from the boundary conditions and Galerkin approach to the
equation (\ref{1}). The approximate solution and the first two derivatives
at the knots can be found from the Eq. (\ref{3}-\ref{4}) as%
\begin{equation}
\begin{tabular}{l}
$U_{i}=U_{N}(x_{i},t)=\alpha _{1}\delta _{i-1}+\delta _{i}+\alpha _{1}\delta
_{i+1},$ \\
$U_{i}^{\prime }=U_{N}^{\prime }(x_{i},t)=\alpha _{2}\delta _{i-1}-\alpha
_{2}\delta _{i+1},$ \\
$U_{i}^{\prime \prime }=U_{N}^{\prime \prime }(x_{i},t)=\alpha _{3}\delta
_{i-1}-2\alpha _{3}\delta _{i}+\alpha _{3}\delta _{i+1}$%
\end{tabular}
\label{5}
\end{equation}%
where $\alpha _{1}=\dfrac{s-ph}{2(phc-s)},\alpha _{2}=\dfrac{p(1-c)}{2(phc-s)%
},\alpha _{3}=\dfrac{p^{2}s}{2(phc-s)}.$

The approximate solution $U_{N}$\ over the element $[x_{m},x_{m+1}]$ can be
written as%
\begin{eqnarray}
U_{N}^{e} &=&\phi _{m-1}\left( x\right) \delta _{m-1}\left( t\right) +\phi
_{m}\left( x\right) \delta _{m}\left( t\right) +\phi _{m+1}\left( x\right)
\delta _{m+1}\left( t\right)  \notag \\
&&+\phi _{m+2}\left( x\right) \delta _{m+2}\left( t\right)  \label{7}
\end{eqnarray}%
where quantities $\delta _{j}\left( t\right) ,j=m-1,...,m+2$ are element
parameters and $\phi _{j}\left( x\right) ,j=m-1,...,m+2$ are known as the
element shape functions.

Over the sample interval $[x_{m},x_{m+1}],$ applying Galerkin approach to
Eq. (\ref{1}) with the test function $\phi _{j}\left( x\right) $ yields the
integral equation:%
\begin{equation}
\underset{x_{m}}{\overset{x_{m+1}}{\int }}\phi _{j}\left( x\right) \left(
u_{t}+uu_{x}-\nu u_{xx}\right) dx.  \label{8}
\end{equation}%
Substitution of the Eq. (\ref{7}) into the integral equation lead to
\begin{eqnarray}
&&\left. \overset{m+2}{\underset{i=m-1}{\sum }}\left( \underset{x_{m}}{%
\overset{x_{m+1}}{\int }}\phi _{j}\phi _{i}dx\right) \overset{\mathbf{%
\bullet }}{\delta }_{i}+\left( \underset{x_{m}}{\overset{x_{m+1}}{\int }}%
\phi _{j}\left( \overset{m+2}{\underset{k=m-1}{\sum }}\delta _{k}\phi
_{k}\right) \phi _{i}^{\prime }dx\right) \delta _{i}\right.  \notag \\
&&\qquad \left. -\nu \left( \underset{x_{m}}{\overset{x_{m+1}}{\int }}\phi
_{j}\phi _{i}^{\prime \prime }dx\right) \delta _{i},\right.  \label{9}
\end{eqnarray}%
where $i,j$ and $k$ take only the values $m-1,m,m+1,m+2$ for $m=0,1,\ldots
,N-1$ and $\overset{\mathbf{\bullet }}{}$ denotes time derivative.

If we denote $A_{ji}^{e},B_{jki}^{e}(\delta ^{e})$ and $C_{ji}^{e}$ by%
\begin{equation}
\begin{tabular}{ll}
$A_{ji}^{e}=\underset{x_{m}}{\overset{x_{m+1}}{\int }}\phi _{j}\phi _{i}dx,$
& $B_{jki}^{e}\left( \delta \right) =\underset{x_{m}}{\overset{x_{m+1}}{%
\int }}\phi _{j}\left( \overset{m+2}{\underset{k=m-1}{\sum }}\delta
_{k}\phi _{k}\right) \phi _{i}^{\prime }dx,$ \\
$C_{ji}^{e}=\underset{x_{m}}{\overset{x_{m+1}}{\int }}\phi _{j}\phi
_{i}^{\prime \prime }dx$ &
\end{tabular}
\label{10}
\end{equation}%
where $\mathbf{A}^{e}$ and $\mathbf{C}^{e}$ are the element matrices of
which dimensions are $4\times 4$ and $\mathbf{B}^{e}\left( \mathbf{\delta }%
^{e}\right) $ is the element matrix with the dimension $4\times 4\times 4$,
Eq.(\ref{9}) can be written in the matrix form as%
\begin{equation}
\mathbf{A}^{e}\overset{\mathbf{\bullet }}{\mathbf{\delta }^{e}}+\left(
\mathbf{B}^{e}\left( \mathbf{\delta }^{e}\right) -\nu \mathbf{C}^{e}\right)
\mathbf{\delta }^{e},  \label{11}
\end{equation}%
where $\mathbf{\delta }^{e}\mathbf{=}\left( \delta _{m-1},...,\delta
_{m+2}\right) ^{T}.$

Gathering the systems (\ref{11}) over all elements, we obtain global system

\begin{equation}
\mathbf{A}\overset{\mathbf{\bullet }}{\mathbf{\delta }}+\left( \mathbf{B}%
\left( \mathbf{\delta }\right) -\nu \mathbf{C}\right) \mathbf{\delta }=0
\label{12}
\end{equation}%
where $\mathbf{A},\mathbf{B}\left( \mathbf{\delta }\right) ,\mathbf{C}$ are
derived from the corresponding element matrices $\mathbf{A}^{e},\mathbf{B}%
^{e}\left( \mathbf{\delta }^{e}\right) ,\mathbf{C}^{e}$, respectively and $%
\mathbf{\delta =}\left( \delta _{-1},...,\delta _{N+1}\right) ^{T}$ contain
all elements parameters.

The unknown parameters $\mathbf{\delta }$ are interpolated between two time
levels $n$ and $n+1$ with the Crank-Nicolson method%
\begin{equation*}
\begin{array}{cc}
\mathbf{\delta }=\dfrac{\mathbf{\delta }^{n+1}+\mathbf{\delta }^{n}}{2}, &
\overset{\mathbf{\bullet }}{\mathbf{\delta }}=\dfrac{\mathbf{\delta }^{n+1}-%
\mathbf{\delta }^{n}}{\Delta t},%
\end{array}%
\end{equation*}%
we obtain iterative formula for the time parameters $\mathbf{\delta }^{n}$%
\begin{equation}
\left[ \mathbf{A+}\frac{\Delta t}{2}\left( \mathbf{B}\left( \mathbf{\delta }%
^{n+1}\right) -\nu \mathbf{C}\right) \right] \mathbf{\delta }^{n+1}=\left[
\mathbf{A-}\frac{\Delta t}{2}\left( \mathbf{B}\left( \mathbf{\delta }%
^{n}\right) -\nu \mathbf{C}\right) \right] \mathbf{\delta }^{n}.  \label{13}
\end{equation}%
The set of equations consist of $\left( N+3\right) $ equations with $\left(
N+3\right) $ unknown parameters. Boundary conditions must be adapted into
the system. Because of the this requirement, initially the first and last
equations are eliminated from the (\ref{13}) and parameters $\delta
_{-1}^{n+1}$ and $\delta _{N+1}^{n+1}$ are substituted in the remaining
system (\ref{13}) by using following equations%
\begin{equation*}
\begin{array}{l}
u\left( a,t\right) =\alpha _{1}\delta _{-1}^{n+1}+\delta _{0}^{n+1}+\alpha
_{1}\delta _{1}^{n+1}=\beta _{1}, \\
u\left( b,t\right) =\alpha _{1}\delta _{N-1}^{n+1}+\delta _{N}^{n+1}+\alpha
_{1}\delta _{N+1}^{n+1}=\beta _{2}%
\end{array}%
\end{equation*}%
which are obtained from the boundary conditions. Thus we obtain a
septa-diagonal matrix with the dimension $\left( N+1\right) \times \left(
N+1\right) $. Since the system (\ref{13}) is an implicit system together
with the nonlinear term $\mathbf{B}\left( \mathbf{\delta }^{n+1}\right) $,
we have used the following inner iteration at each time step $(n+1)\Delta t$
to work up solutions:%
\begin{equation}
(\mathbf{\delta }^{\ast }\mathbf{)}^{n+1}=\mathbf{\delta }^{n}+\dfrac{(%
\mathbf{\delta }^{n}-\mathbf{\delta }^{n-1})}{2}.  \label{14}
\end{equation}%
We use the above iteration three times to find the new approximation $(%
\mathbf{\delta }^{\ast }\mathbf{)}^{n+1}$ for the parameters $\mathbf{\delta
}^{n+1}$ to recover solutions at time step $(n+1)\Delta t$.

To start evolution of the iterative system for the unknown $\mathbf{\delta }%
^{n}$, the vector of initial parameters $\mathbf{\delta }^{0}$must be
determined by using the following initial and boundary conditions:%
\begin{equation}
\begin{tabular}{l}
$u_{0}^{\prime }(x_{0},0)=\dfrac{p\left( 1-c\right) }{2\left( phc-s\right) }%
\delta _{-1}+\dfrac{p\left( c-1\right) }{2\left( phc-s\right) }\delta _{1}$
\\
$u\left( x_{m},0\right) =\dfrac{s-ph}{2\left( phc-s\right) }\delta
_{m-1}+\delta _{m}+\dfrac{s-ph}{2\left( phc-s\right) }\delta _{m+1},$ $%
m=0,...,N.$ \\
$u^{\prime }\left( x_{N},0\right) =\dfrac{p\left( 1-c\right) }{2\left(
phc-s\right) }\delta _{N-1}+\dfrac{p\left( c-1\right) }{2\left( phc-s\right)
}\delta _{N+1}$%
\end{tabular}
\label{15}
\end{equation}%
The solution of matrix equation (\ref{15}) with the dimensions $\left(
N+1\right) \times \left( N+1\right) $ is obtained by the way of Thomas
algorithm. Once $\mathbf{\delta }^{0}$ is determined, we can start the
iteration of the system to find the parameters $\mathbf{\delta }^{n}$ at
time $t^{n}=n\Delta t.$ Approximate solutions at the knots is found from the
Eq.(\ref{5}) and solution over the intervals $[x_{m},x_{m+1}]$ is determined
from the Eq.(\ref{7}).

\section{Test Problems}

The robustness of the algorithm is shown by studying two test problems.
Error is measured by the maximum error norm;%
\begin{equation}
L_{\infty }=\left \Vert u^{\text{exact}}-u^{\text{numeric}}\right \Vert
_{\infty }=\max_{0\leq j\leq N}\left \vert u_{j}^{\text{exact}}-u_{j}^{\text{%
numeric}}\right \vert .  \label{16}
\end{equation}%
The free parameter $p$ of the exponential B-spline is found by scanning the
predetermined interval with very small increment.

\textbf{(a)}\ A shock propagation solution of the Burgers' equation is%
\begin{equation}
u(x,t)=\dfrac{x/t}{1+\sqrt{t/t_{0}}\exp (x^{2}/(4\nu t))},\text{\quad }t\geq
1,  \label{17}
\end{equation}%
where $t_{0}=\exp (1/(8\nu ))$.\ The sharpness of the solutions increases
with selection of the smaller $\nu $.

\begin{figure}[ht]
\centering\includegraphics[scale=0.5]{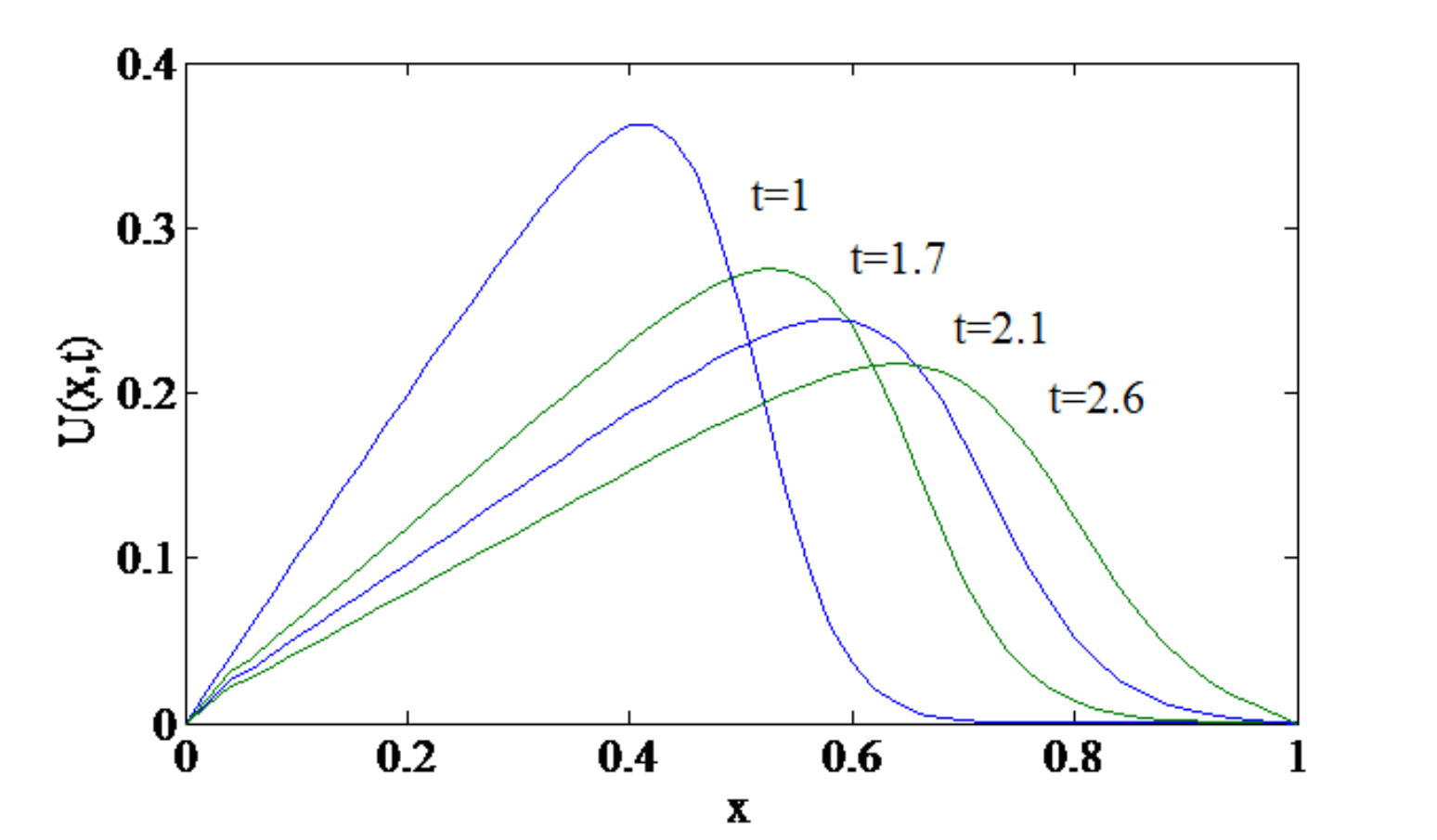}
\caption{{}{\protect \ Solutions for $\nu =0.01\text{, }h=0.02$, $p=0.005111$}}
\label{fig1}
\end{figure}

\begin{figure}[ht]
\centering\includegraphics[scale=0.5]{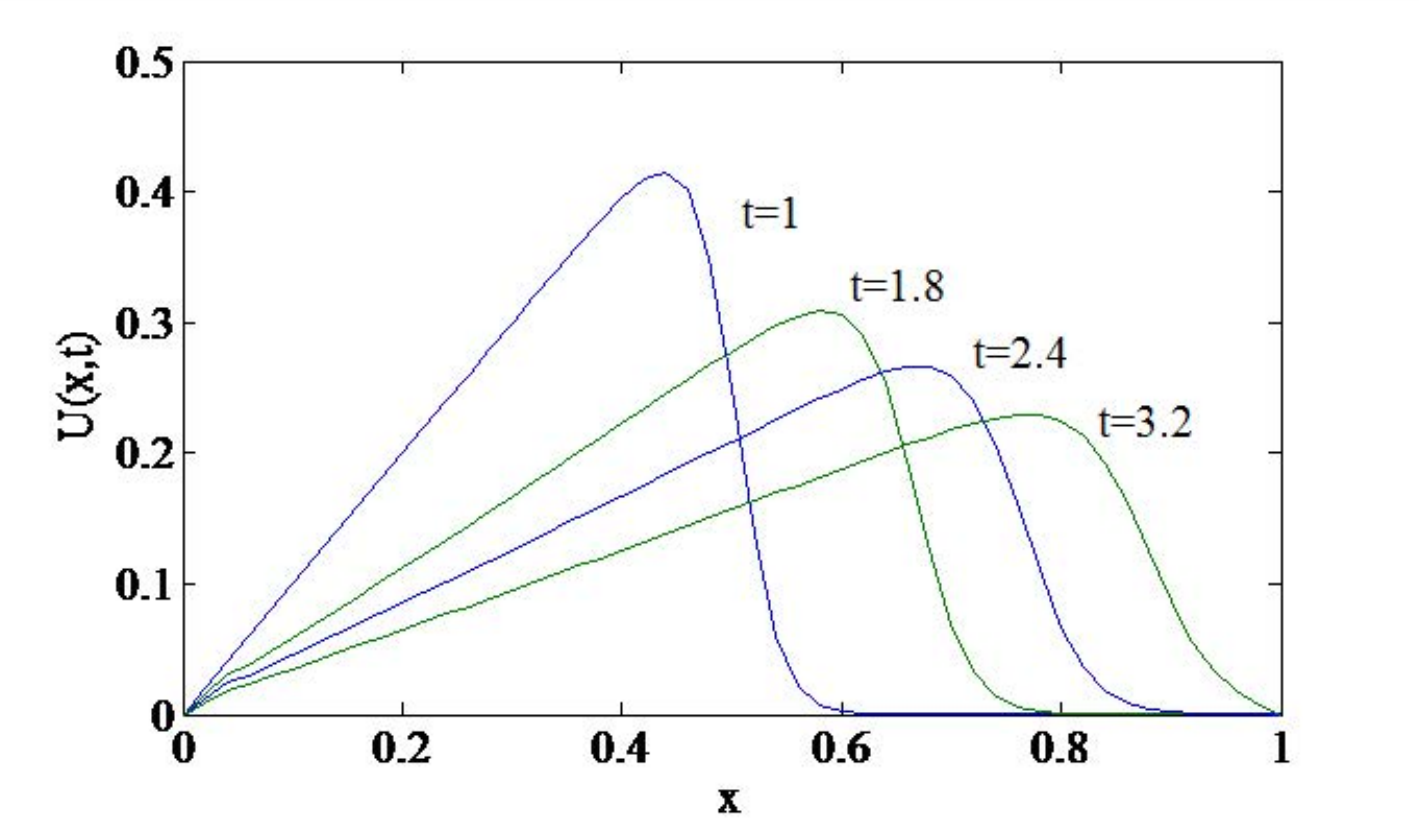}
\caption{{}{\protect \ Solutions for $\nu =0.005\text{, }h=0.02$, $p=0.000739$}}
\label{fig2}
\end{figure}

\begin{figure}[ht]
\centering\includegraphics[scale=0.5]{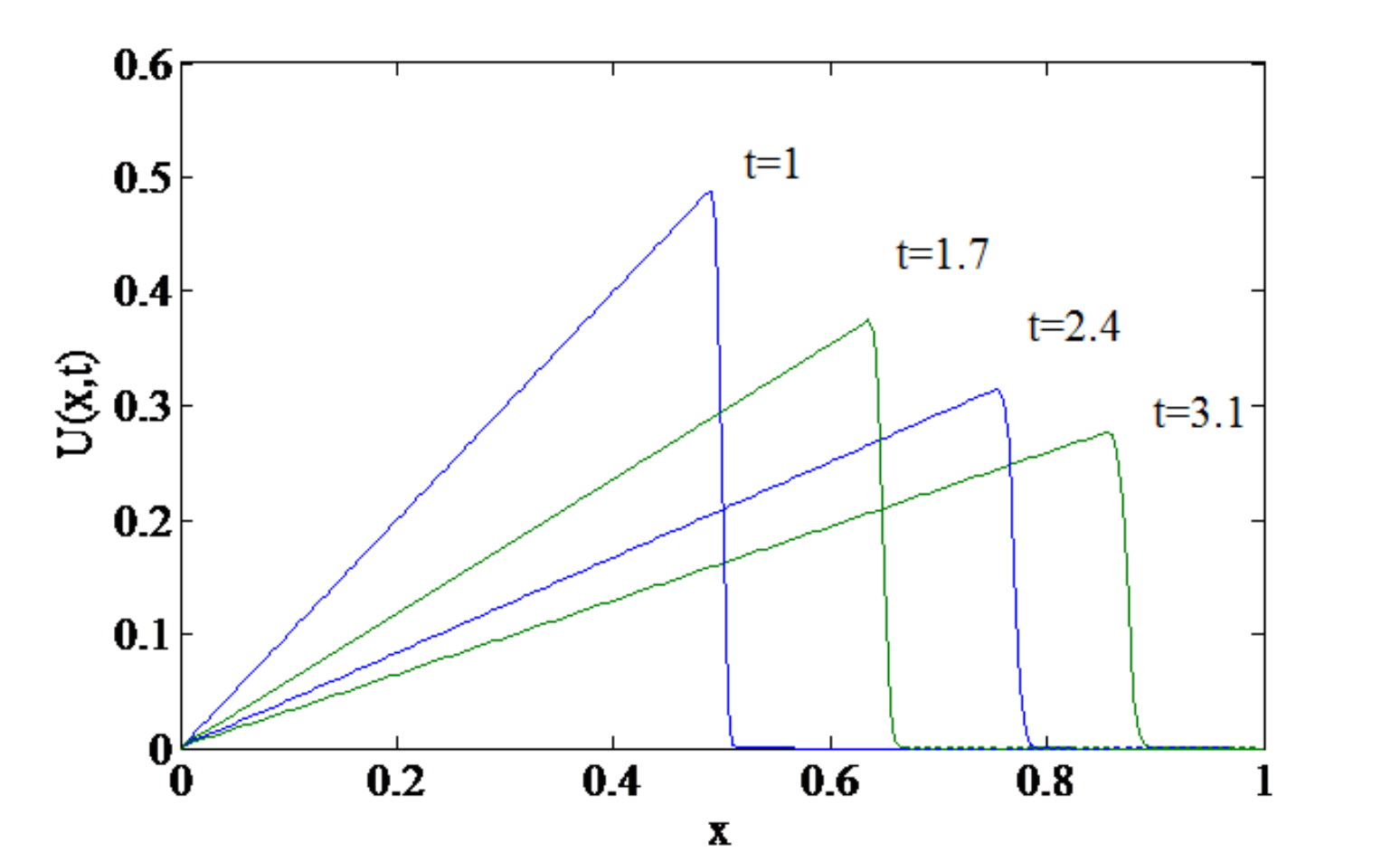}
\caption{{}{\protect \ Solutions for $\nu=0.0005\text{, }h=0.005$, $p=0.005941$}}
\label{fig3}
\end{figure}

Substitution of the $t=1$ in Eq. (\ref{17}) gives the initial condition. The
boundary conditions $u(0,t)=0$ and $u(1,t)=0$ are used. Computations are
performed with parameters $\nu =0.0005,$ $0.005,$ $0.01,$ $h=0.02,$ $0.005$
and $\Delta t=0.01$ over the solution domain $[0,1]$. As time increases,
shock evaluation is observed and some graphical solutions are drawn in 
Figs. \ref{fig1}-\ref{fig3} for various viscosity values and space steps. For $\nu =0.01,$ algorithm
produces smoother shock during run time. With decreasing values of $\nu $,
as seen in the Figs. \ref{fig1}-\ref{fig3} the steepening occurs. For the smaller viscosity
constant $\nu =0.0005$ the sharper shock is observed and steepness of
numerical solution is kept almost unchanged during the program run. The
results obtained by present scheme can be compared with ones given in the
works \cite{bs1,ra,bs2,bs3,bs4} through the computation of\ $L_{\infty }$
error norm at various times in the Table 2.%

\begin{equation*}
\begin{tabular}{llll}
\multicolumn{4}{l}{$\text{Table 2: Comparison of numerical results at
different times}$} \\ \hline \hline
& $L_{\infty }\times 10^{3}$ & $L_{\infty }\times 10^{3}$ & $L_{\infty
}\times 10^{3}$ \\ \cline{2-4}
$h=0.005\text{, }\nu =0.005$ & $t=1.7$ & $t=2.4$ & $t=3.1$ \\ \hline
$\text{Present (}p=0.005941\text{)}$ & $3.15776$ & $2.33757$ & $4.79061$ \\
$\text{Ref.\cite{bs1} (QBGM)}$ & $1.20755$ & $0.80187$ & $4.79061$ \\
$\text{Ref.\cite{bs2} (QBCM1)}$ & $0.06192$ & $0.05882$ & $4.43469$ \\
$\text{Ref.\cite{bs3} (QBCA1)}$ & $1.21175$ & $0.80771$ & $4.79061$ \\
$\text{Ref.\cite{bs4}}$ & $0.04284$ & $0.06464$ & $4.79061$ \\
&  &  &  \\
$h=0.02\text{, }\nu =0.005$ & $t=1.8$ & $t=2.4$ & $t=3.2$ \\ \hline
$\text{Present (}p=0.000739\text{)}$ & $8.26075$ & $7.42050$ & $7.49146$ \\
$\text{Ref.\cite{ra}}$ & $2.47189$ & $2.16784$ & $7.49146$ \\
$\text{Ref.\cite{bs2} (QBCM1)}$ & $0.54058$ & $0.39241$ & $5.54899$ \\
$\text{Ref.\cite{bs3} (QBCA1)}$ & $1.15263$ & $0.80008$ & $7.49147$ \\
$\text{Ref.\cite{bs4}}$ & $0.03546$ & $0.06464$ & $7.49147$ \\
&  &  &  \\
$h=0.02\text{, }\nu =0.01$ & $t=1.7$ & $t=2.1$ & $t=2.6$ \\ \hline
$\text{Present (}p=0.005111\text{)}$ & $8.08651$ & $7.53518$ & $8.06798$ \\
$\text{Ref.\cite{ra}}$ & $3.13476$ & $2.66986$ & $8.06798$ \\
$\text{Ref.\cite{bs2} (QBCM1)}$ & $0.40431$ & $0.86363$ & $6.69425$ \\
$\text{Ref.\cite{bs3} (QBCA1)}$ & $0.47456$ & $1.14759$ & $8.06798$ \\
$\text{Ref.\cite{bs4}}$ & $0.09592$ & $1.14760$ & $8.06799$ \\ \hline \hline
\end{tabular}%
\end{equation*}

\begin{figure}[ht]
\centering\includegraphics[scale=0.5]{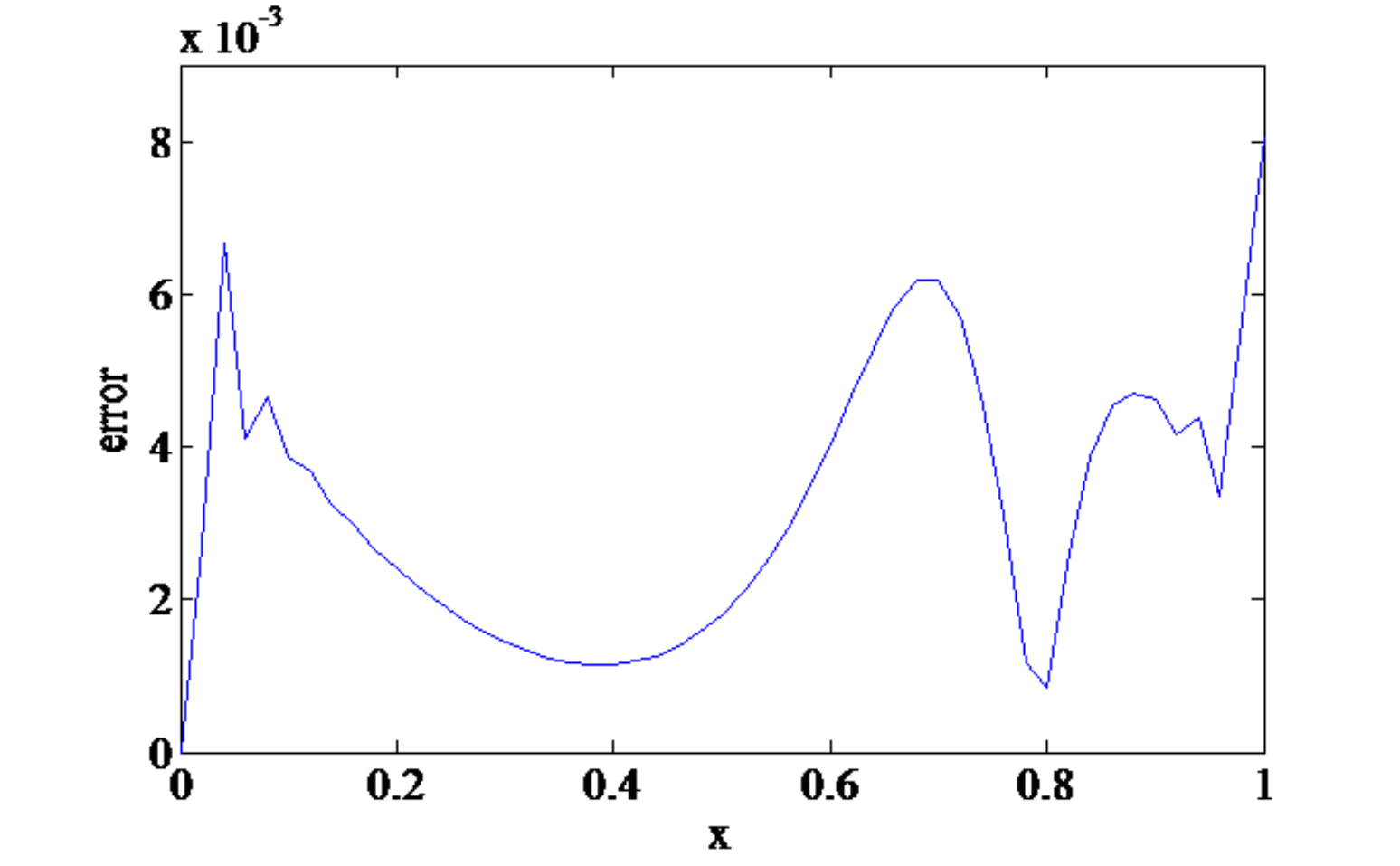}
\caption{{}{\protect \ Absolute error for $\nu =0.01\text{, }h=0.02\text{, }p=0.005111$}}
\label{fig4}
\end{figure}

\begin{figure}[ht]
\centering\includegraphics[scale=0.5]{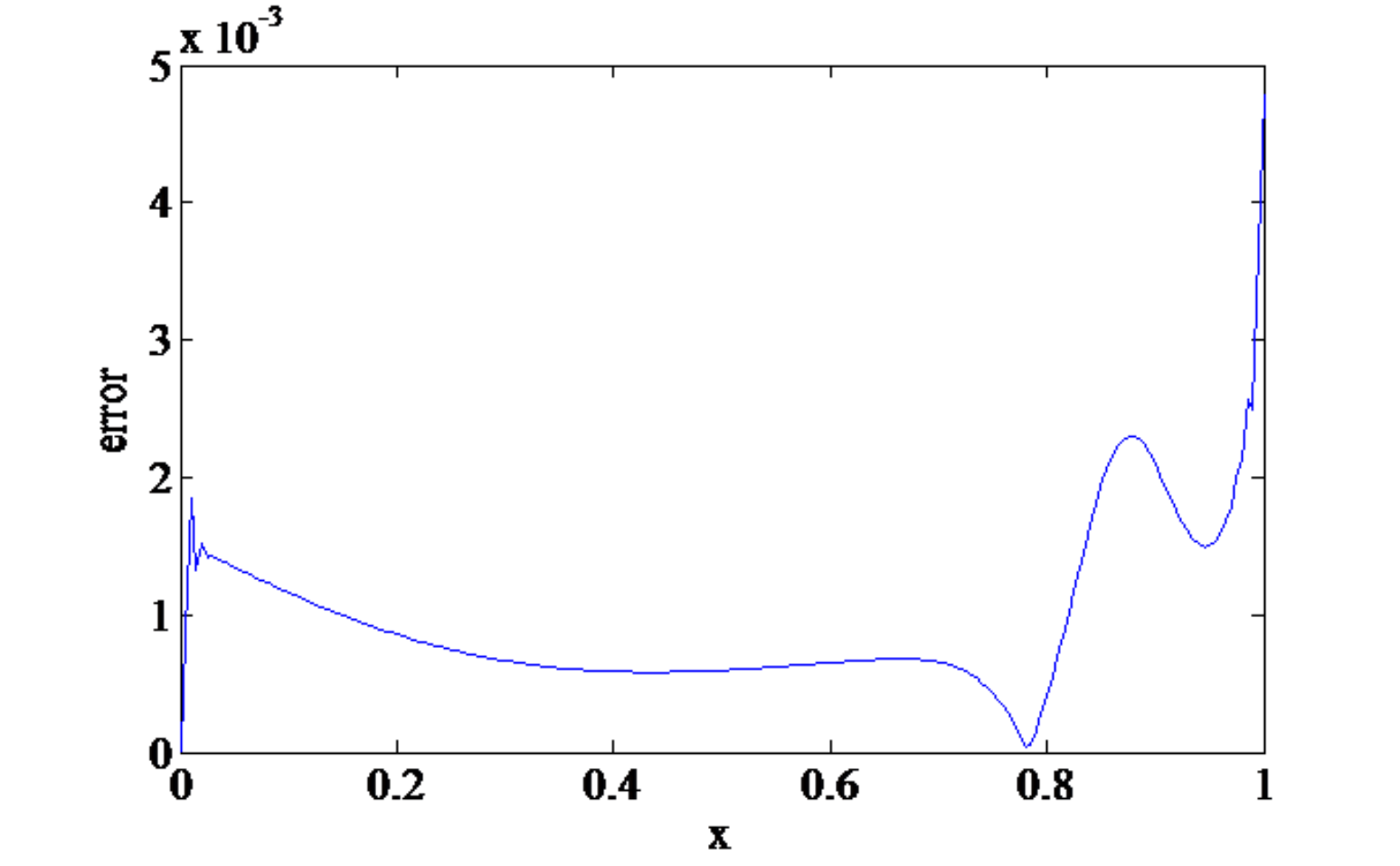}
\caption{{}{\protect \ Absolute error for $\nu =0.005\text{, }h=0.005\text{, }p=0.005941$}}
\label{fig5}
\end{figure}

\begin{figure}[ht]
\centering\includegraphics[scale=0.5]{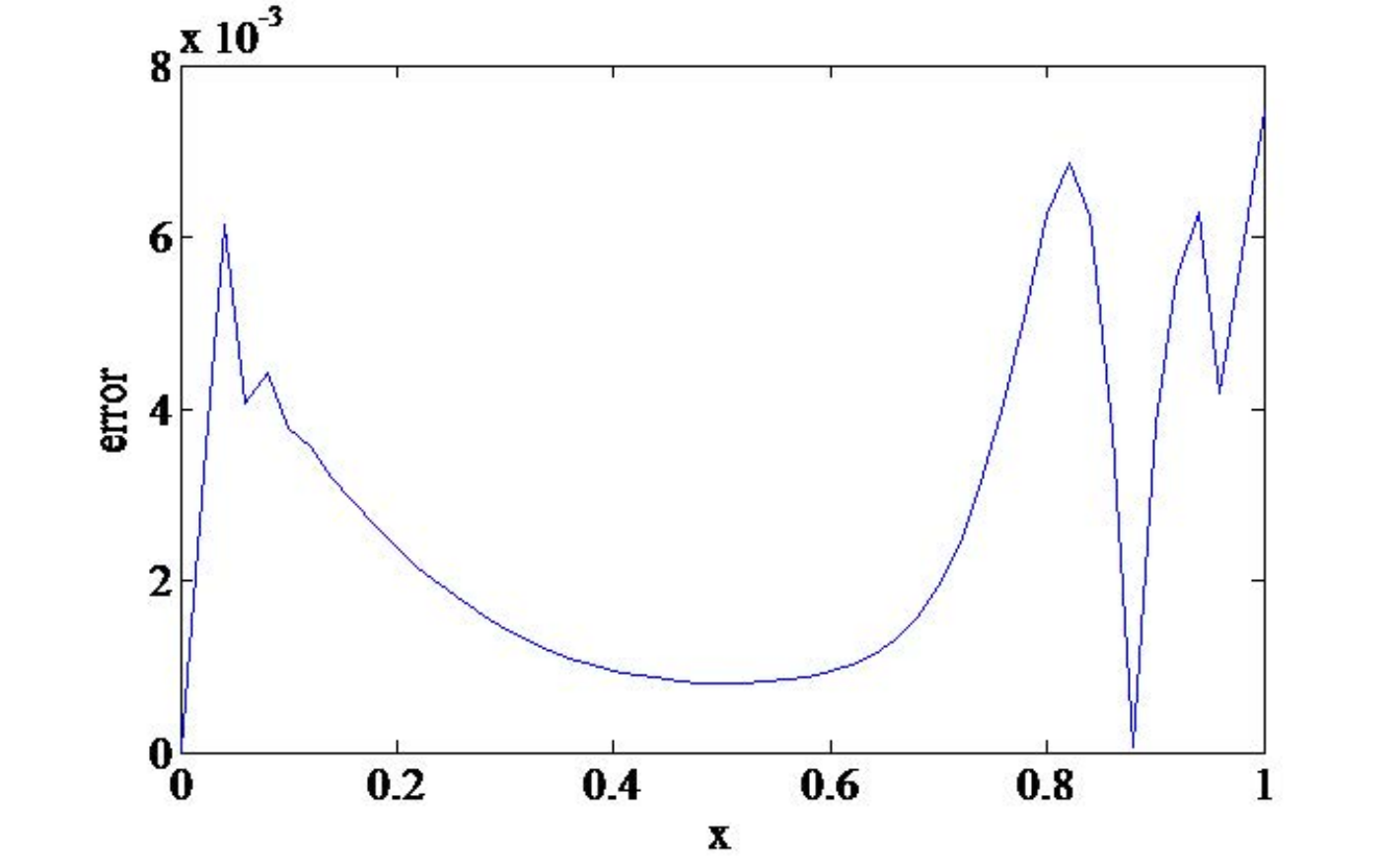}
\caption{{}{\protect \ Absolute error for $\nu =0.005\text{, }h=0.02\text{, }p=0.000739$}}
\label{fig6}
\end{figure}

\begin{figure}[ht]
\centering\includegraphics[scale=0.5]{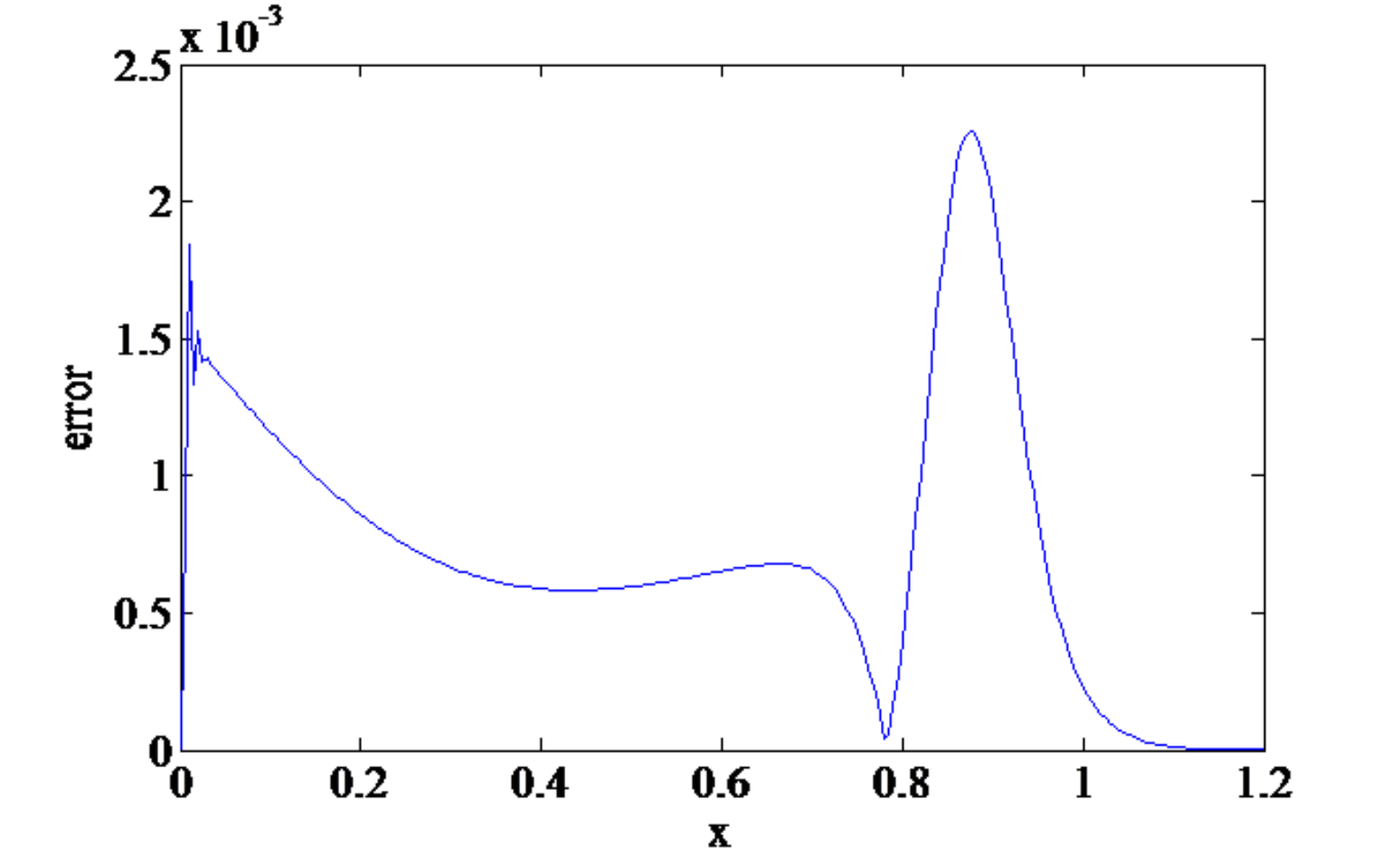}
\caption{{}{\protect \ Absolute error for $\nu =0.005\text{, }h=0.005\text{, }$%
$p=0.005941\text{ over }[0,1.2]$}}
\label{fig7}
\end{figure}

The absolute error distributions between the analytical and numerical
solutions are drawn in Figs. \ref{fig4}-\ref{fig6} for various viscosity values and space
steps. In these figures, the highest error appears about the right hand
boundary position. We run the program again over the extended domain $\left[
0,1.2\right] $ with parameters $\nu =0.005,$ $h=0.005$, the highest error is
reduced at the right boundary seen in the Fig. \ref{fig7} and $L_{\infty }$ error
norm decrease from $4.790609\times 10^{-3}$ to $2.259598\times 10^{-3}$ at
time $t=3.1$.

\textbf{(b)} A well-known analytical solution of Burgers' equation is%
\begin{equation}
u(x,t)=\dfrac{\alpha +\mu +\left( \mu -\alpha \right) \exp \eta }{1+\exp
\eta },\text{\quad }0\leq x\leq 1,\text{ }t\geq 0,  \label{18}
\end{equation}%
where $\eta =\dfrac{\alpha \left( x-\mu t-\gamma \right) }{\nu }$. $\alpha ,$
$\mu $ and $\gamma $ are constants. Parameters \linebreak $\alpha =0.4,$ $%
\mu =0.6$ and $\gamma =0.125$ are used to coincide with the some previous
studies. This solution involves a travelling wave and move to the right with
speed $\mu $. Initial condition is obtained from Eq. (\ref{18}) when $t=0$.
The boundary conditions are $u(0,t)=1,$ $u(1,t)=0.2$ for $t>0$.

The calculation is performed with time step $\Delta t=0.01$, space step $%
h=1/36$ and viscosity coefficient $\nu =0.01$. The program is run up to time
$t=0.5$. We have found $L_{\infty }=6.73543978\times 10^{-4}$ for the
exponential B-spline Galerkin method at time $t=0.5$ documented in Table 3
with results of the quadratic B-spline Galerkin method \cite{bs1}, the
quartic B-spline collocation method \cite{bs2}, the quintic B-spline
collocation method \cite{bs3} and the quartic B-spline Galerkin method \cite%
{bs4}.

\begin{figure}[ht]
\centering\includegraphics[scale=0.5]{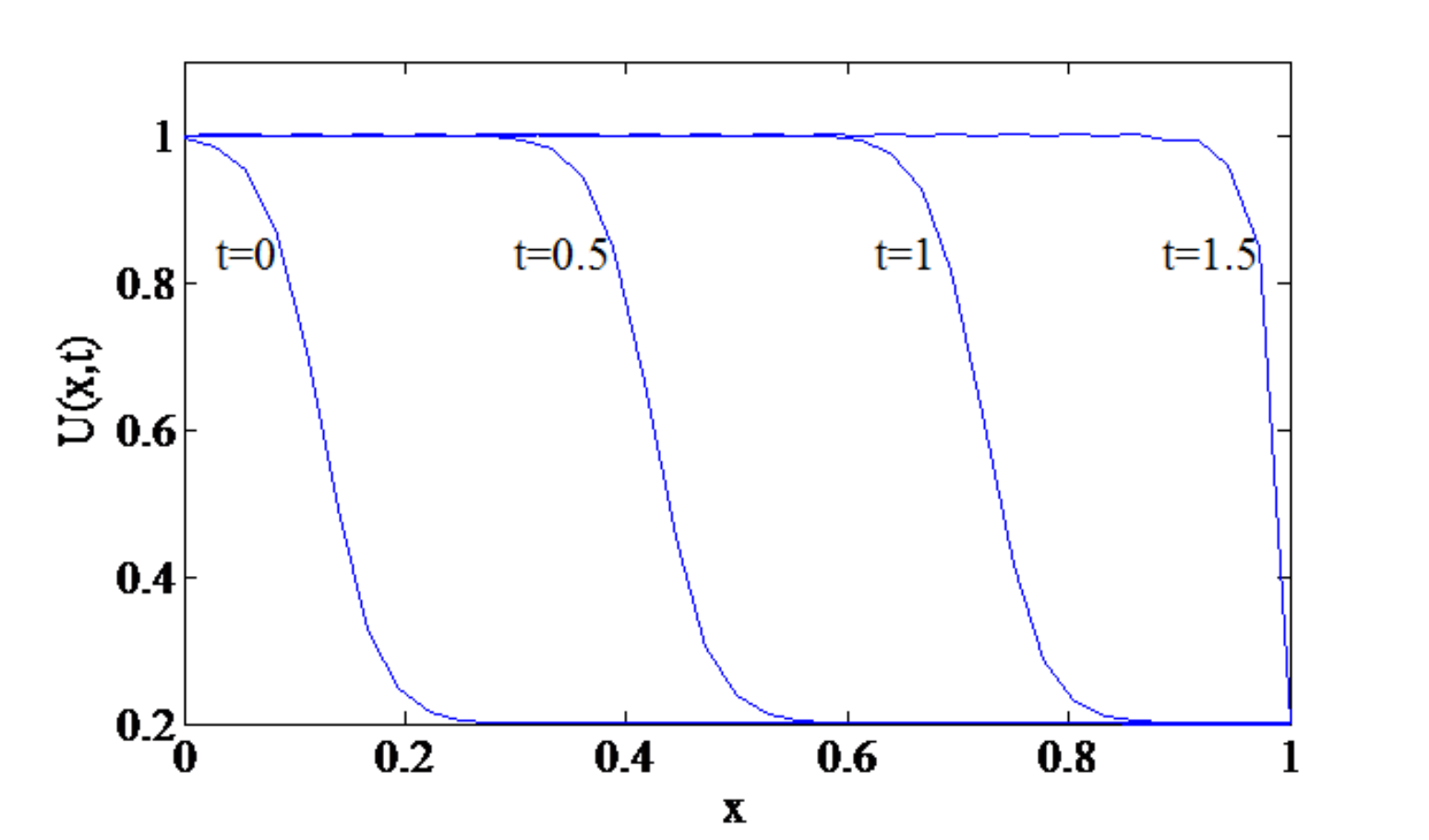}
\caption{{}{\protect \ Solutions for $\nu =0.01$}}
\label{fig8}
\end{figure}

\begin{figure}[ht]
\centering\includegraphics[scale=0.5]{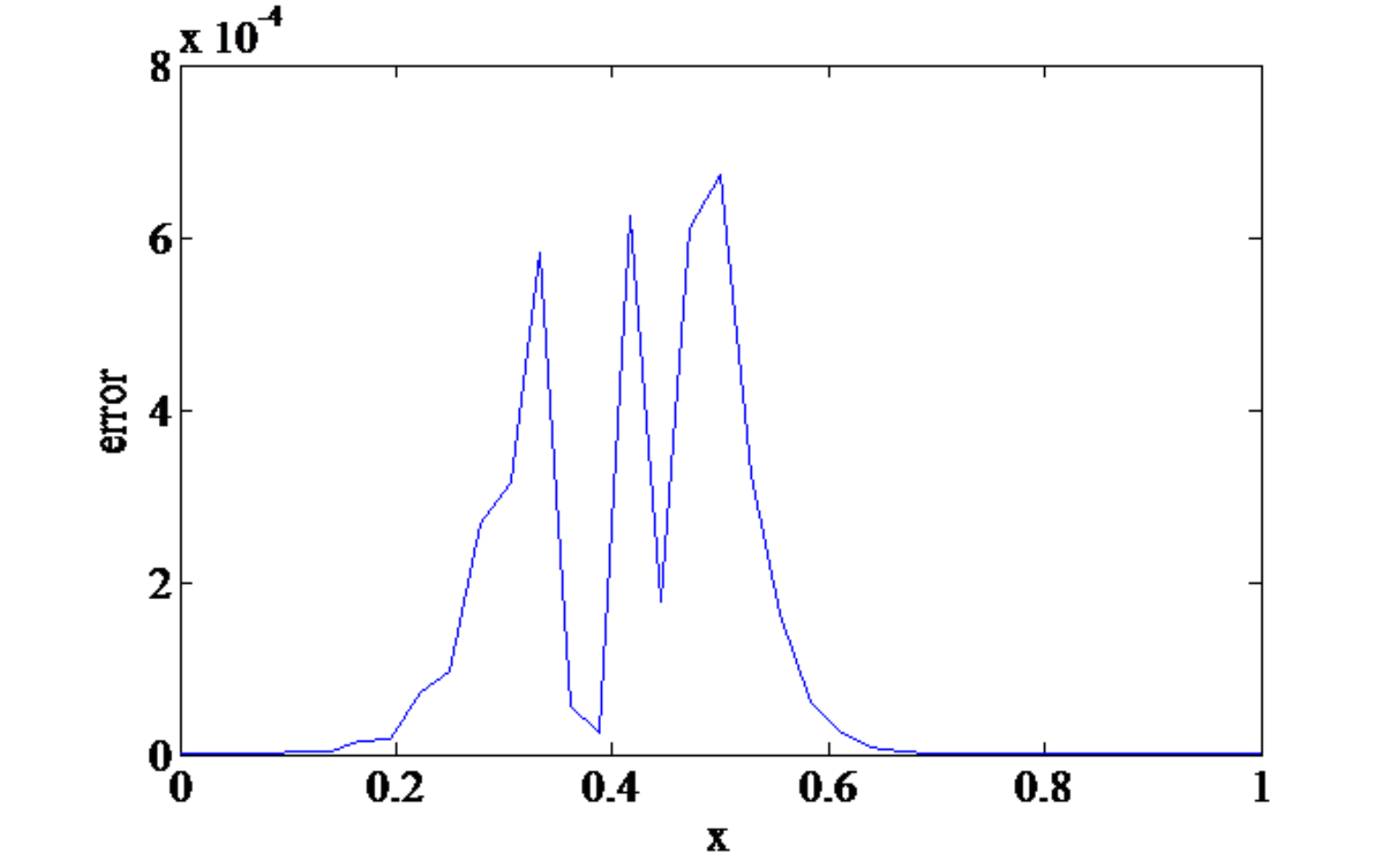}
\caption{{}{\protect \ Absolute error for $\nu =0.01$}}
\label{fig9}
\end{figure}

The numerical solution obtained by the presented schemes gives better
results than the others. The profiles of initial wave and solution at some
times are depicted in Fig. \ref{fig8}. Error variations of the schemes are given in
Fig. \ref{fig9} at time $t=0.5$.%

\begin{equation*}
\begin{tabular}{lll}
\multicolumn{3}{l}{Table 3: Comparison of results at $t=0.5$ for $h=1/36$, $%
\nu =0.01$} \\ \hline \hline
$\text{Method}$ &  & $L_{{\footnotesize \infty }}\times 10^{{\footnotesize 3}%
}$ \\ \hline
$\text{Present (}p=0.002323$) &  & $0.67354$ \\
$\text{Ref.\cite{bs1} (QBGM)}$ &  & $6.35489$ \\
$\text{Ref.\cite{bs2} (QBCM1)}$ &  & $3.03817$ \\
$\text{Ref.\cite{bs3} (QBCA1)}$ &  & $5.78454$ \\
$\text{Ref.\cite{bs4} (QBGM)}$ &  & $1.44$ \\ \hline \hline
\end{tabular}%
\end{equation*}

\section{Conclusion}

In this paper, we investigate the utility of the exponential B-spline in the
Galerkin algorithm for solving the Burgers' equation. The efficiency of the
method is tested for a shock propagation solution and a travelling solution
of the Burgers' equation. For the first test problem, solutions found with
the present methods are in good agreement with the results obtained by
previous studies. In the second test problem, present method leads to
accurate results than all of the others. In conclude, the numerical
algorithm in which the exponential B-spline functions are used, performs
well compared with other existing numerical methods for the solution of
Burgers' equation.\bigskip

\noindent \textbf{Acknowledgements:} The authors are grateful to The
Scientific and Technological Research Council of Turkey for financial
support given for the project 113F394.

\end{document}